\documentclass[a4paper,11pt]{article}

\usepackage{amsmath}
\usepackage{amsthm}
\usepackage[all]{xy}
\usepackage{amssymb}
\usepackage{mathrsfs}
\usepackage{enumerate}
\usepackage[usenames,dvipsnames]{color}

\usepackage{fullpage}
\usepackage{amsmath,amsfonts,amssymb}
\makeatletter
\newcommand{\thorn}{{\fontencoding{T1}\selectfont\th}}
\newcommand{\@indepsymbol}[2]{#1\setbox0=\hbox{$#1x$}\kern\wd0\hbox to 0pt{\hss$#1\mid$\hss}\lower.9\ht0\hbox to 0pt{\hss$#1\smile$\hss}\kern\wd0}
\newcommand{\@nindepsymbol}[2]{#1\setbox0=\hbox{$#1x$}\kern\wd0\hbox to 0pt{\mathchardef
	\nn=12854\hss$#1\nn$\kern1.4\wd0\hss}\hbox to
	0pt{\hss$#1\mid$\hss}\lower.9\ht0 \hbox to
	0pt{\hss$#1\smile$\hss}\kern\wd0}
\newcommand{\ind}[1][]{\mathop{\mathpalette\@indepsymbol{}^{\!\!\!\!\rlap{$\scriptstyle\textnormal{#1}$}\,\,\,\,}}}
\newcommand{\nind}[1][]{\mathop{\mathpalette\@nindepsymbol{}^{\!\!\!\rlap{$\scriptstyle\textnormal{#1}$}\,\,\,}}}
\newcommand{\@Ind}[1][]{\mathpalette\@indepsymbol{}^{\!\!\!\!\mbox{$\scriptstyle\textnormal{#1}$}}}
\newcommand{\Ind}[1][]{\@Ind[\ \,]}
\newcommand{\newind}[4]{
	\newcommand{#1}{{\!\@Ind[#4]}}
	\newcommand{#2}{\ind[#4]}
	\newcommand{#3}{\nind[#4]}
}
\newind{\thInd}{\thind}{\nthind}{\thorn}
\renewcommand{\thInd}{\text{$\@Ind[\thorn]$\;}}
\makeatother

\newtheorem{thm}{Theorem}[section]

\newtheorem{conj}[thm]{Conjecture}

\theoremstyle{definition}

\theoremstyle{remark}

\theoremstyle{remark}

\theoremstyle{remark}

\title{A special case of quasiminimality}

\author{Gareth Boxall}

\date{\today}

\begin{document}

\maketitle

\begin{abstract} We deduce a special case of Zilber's quasiminimality conjecture, for the complex exponential field, from work of Henson and Rubel. Specifically, we deal with those subsets of $\mathbb{C}$ defined by formulas of the form $\exists\bar{y}(P(x,\bar{y})=0)$, where $P$ is a term formed from the language $\{+,\times,\exp\}$ together with parameters from $\mathbb{C}$. 
\end{abstract}

\section{Introduction}

Boris Zilber's quasiminimality conjecture is the following.

\begin{conj}\label{QM}
If $X\subseteq\mathbb{C}$ is definable in the structure $(\mathbb{C},+,\times,\exp)$ then either $X$ or $\mathbb{C}\setminus X$ is countable.
\end{conj} 

Here $\exp:\mathbb{C}\rightarrow\mathbb{C}$ is the function which sends $z$ to $e^z$. A structure is said to be minimal if every definable subset of its underlying set is finite or cofinite. For example, $(\mathbb{C},+,\times)$ is minimal by quantifier elimination and basic facts of polynomials. It is well known that $(\mathbb{C},+,\times,\exp)$ is not minimal because $\mathbb{Q}$ is definable in it. The weakening of minimality in which ``finite" is replaced by ``countable" is called quasiminimality. 

Conjecture \ref{QM} remains open. However there has been significant progress in proving it conditionally. In \cite{BK} Bays and Kirby establish it under the assumption that $\mathbb{C}$ is exponentially-algebraically closed. Wilkie has outlined a different approach via a conjecture of his that certain implicitly defined functions will analytically continue to all but countably many points of $\mathbb{C}$. This is discussed, for example, in \cite{W} where the analytic continuation conjecture is stated precisely in a simple case and it is mentioned that a suitable generalisation would imply all of Conjecture 1.1.

Despite this progress, not many special cases of Conjecture 1.1 are known unconditionally. It is trivial when $X$ is defined by a quantifier-free formula, since then $X$ will be a boolean combination of zero-sets of entire functions and such sets are either countable or $\mathbb{C}$. But introducing quantifiers causes problems. The case where $X$ is defined by a formula consisting of some universal quantification followed by an atomic formula is easy and presumably known ($\forall\bar{y}(P(x,\bar{y})=0)$ is equivalent to $(\forall \bar{y}\in\mathbb{Q}^m)(P(x,\bar{y})=0)$ and we are done). In this note we establish the case where $X$ is defined by a formula consisting of some existential quantification followed by an atomic formula. The proof uses work of Henson and Rubel which employs Nevanlinna theory in the study of exponential polynomials.

For all $A\subseteq \mathbb{C}$, let $L_{A}$ be the language $\{+,\times,\exp\}$ together with a constant symbol for each element of $A$. Let $\Sigma$ be the set of all $L_{\mathbb{C}}$-terms. Let $t=t(x_1,...,x_n)\in\Sigma$. It is shown in \cite{HR} that if the function from $\mathbb{C}^n$ to $\mathbb{C}$ defined by $t$ has no zeroes then there is some $s\in\Sigma$ such that $t=\exp(s)$, where here, and throughout, equality between terms means equality of the complex functions they define. For each $A\subseteq\mathbb{C}$, let $\Sigma_A\subseteq \Sigma$ consist of those terms whose constant symbols are in $L_A$. We would like Henson and Rubel's result with $\Sigma_A$ in place of $\Sigma$. It turns out this is possible (using their proof) up to a constant multiple. 

\begin{thm}\label{HR}
Let $A\subseteq\mathbb{C}$ and let $t=t(x_1,...,x_n)\in \Sigma_A$. Assume the function from $\mathbb{C}^n$ to $\mathbb{C}$ defined by $t$ has no zeroes. Then there exist $s\in\Sigma_A$ and a constant symbol $c\in L_{\mathbb{C}}$ such that $t=c(\exp(s))$. 
\end{thm}

Using this we prove the following.

\begin{thm}\label{sc}
Let $P(x,\bar{y})\in\Sigma$, where $x$ is a single variable and $\bar{y}$ is a tuple of variables. Let $X$ be the subset of $\mathbb{C}$ defined by $\exists\bar{y}(P(x,\bar{y})=0)$. Then either $X$ or $\mathbb{C}\setminus X$ is countable.  
\end{thm}

We actually prove that the formula $\exists\bar{y}(P(x,\bar{y})=0)$ is equivalent to a countable boolean combination of formulas of the form  $(\exists \bar{y}\in \mathbb{Q}^m)\varphi(x,\bar{y})$, where $\varphi$ is a quantifier-free $L_\mathbb{C}$-formula which has no new parameters. Our argument works just as well for a tuple $\bar{x}$ in place of the single variable $x$.

\begin{thm}\label{final}
Let $A\subseteq\mathbb{C}$. Every formula of the form $\exists\bar{y}(P(\bar{x},\bar{y})=0)$, where $P(\bar{x},\bar{y})\in \Sigma_A$, is equivalent in $(\mathbb{C},+,\times,\exp)$ to a countable boolean combination of formulas of the form $(\exists \bar{y}\in \mathbb{Q}^m)\varphi(\bar{x},\bar{y})$, where $\varphi$ is a quantifier-free $L_A$-formula. 
\end{thm}

\section{Work of Hensen and Rubel}

One proves Theorem \ref{HR} by inspecting the proof of Theorem 5.4 in \cite{HR} and noticing that the argument given there establishes the result. Specifically, in the argument which we now present, we differ from \cite{HR} only in requiring $s_1,...,s_k\in \Sigma_A$ (as opposed to $s_1,...,s_k\in\Sigma$) and our only contribution is to check that reasoning within the proof which allows them to lower their $k$ allows us to lower our one too. 

Let $p_1,...,p_k\in \Sigma$ and $s_1,...,s_k\in \Sigma_A$ such that $t=p_1\exp(s_1)+...+p_k\exp(s_k)$ and $\exp$ does not occur in $p_1,...,p_k$ (so the $p_i$'s correspond to polynomials). Assume $k$ is minimal. If $k=1$ then $p_1$ defines a polynomial function with no zeroes which is therefore constant and we are done. It remains to show that indeed $k=1$. 

Suppose $k>1$. The functions $p_1\exp(s_1),...,p_k\exp(s_k)$ are linearly independent over $\mathbb{C}$, since otherwise we could drop one of them from the sum and adjust the remaining terms by constant factors so as to again obtain $t$ and thereby contradict the minimality of $k$. 

The proof proceeds, just as in \cite{HR}, to the conclusion that $\frac{\exp(s_i)}{\exp(s_j)}$ is constant for some distinct $i,j\in\{1,...,m\}$. This contradicts the minimality of $k$, since the expression $p_i\exp(s_i)+p_j\exp(s_j)$ may then be simplified to $(p_i+cp_j)\exp(s_i)$ for some constant $c$.

\section{Proofs}

We shall not be too strict about the distinction between elements of $\Sigma$ and the functions or complex numbers they define or represent. Let $A\subseteq\mathbb{C}$ and $P(\bar{x},\bar{y})\in\Sigma_A$. Say $\bar{x}=(x_1,...,x_n)$ and $\bar{y}=(y_1,...,y_m)$. For each $Q(\bar{x},\bar{y})\in\Sigma_A$ and $i\in\{1,...,m\}$, let $$Q_i(\bar{x},\bar{y})=\frac{\partial}{\partial y_i}\left(\frac{P(\bar{x},\bar{y})}{\exp(Q(\bar{x},\bar{y}))}\right).$$ For each $\bar{a}\in\mathbb{C}^n$, the following are equivalent (the equivalence of the first two coming from Theorem \ref{HR}).

\begin{enumerate}

\item $P(\bar{a},\bar{y})$ has no zeroes.

\item $P(\bar{a},\bar{y})=c(\exp(Q(\bar{a},\bar{y})))$ for some $Q(\bar{x},\bar{y})\in \Sigma_A$ and non-zero $c\in\mathbb{C}$. 

\item $P(\bar{a},\bar{y})$ is not constantly zero and there exists $Q(\bar{x},\bar{y})\in\Sigma_A$ such that, for each $i\in\{1,...,m\}$, $Q_i(\bar{a},\bar{y})$ is constantly zero.

\end{enumerate}

The third of these conditions is equivalent to the existence of $Q(\bar{x},\bar{y})\in\Sigma_A$ such that

\begin{enumerate}

\item[i.] $P(\bar{a},\bar{q})\neq 0$ for some $\bar{q}\in\mathbb{Q}^m$ and

\item[ii.] $Q_i(\bar{a},\bar{q})=0$ for all $\bar{q}\in\mathbb{Q}^m$ and $i\in\{1,...,m\}$.

\end{enumerate}

For all $Q(\bar{x},\bar{y})\in\Sigma_A$ and $i\in\{1,...,m\}$,  $Q_i(\bar{x},\bar{y})$ is a fraction of members of $\Sigma_A$. Since $\Sigma_A$ is countable, we obtain Theorem \ref{final}. Given that $\mathbb{Q}$ is countable and Conjecture \ref{QM} holds for quantifier-free definable sets, this implies Theorem \ref{sc}.

\bibliographystyle{abbrv}

\begin{thebibliography}{99}

\bibitem{BK} Bays, M. and Kirby, J., Pseudo-exponential maps, variants, and quasiminimality, arXiv:1512.04262v2, 2016.

\bibitem{HR} Henson, C. W. and Rubel, L. A., Some applications of Nevanlinna theory to mathematical logic: Identities of exponential functions, {\em Trans. Amer. Math. Soc.}, 282: 1--32, 1984.

\bibitem{W} Wilkie, A., Some results and problems on complex germs with definable Mittag-Leffler stars, {\em Notre Dame J. Form. Log.}, 54: 603--610, 2013. 








\end{thebibliography}

\end{document}